\documentclass[10pt,a4paper]{article}
\usepackage[utf8]{inputenc}
\usepackage{amsmath}
\usepackage{amsfonts}
\usepackage{amssymb}
\usepackage{graphicx}
\usepackage{float}

\usepackage{amsthm}
\usepackage{bbm}
\usepackage{xcolor}
\usepackage{mathtools}

\newtheorem{theorem}{Theorem}[section]
\newtheorem{lemma}[theorem]{Lemma}
\newtheorem{proposition}[theorem]{Proposition}

\theoremstyle{definition}
\newtheorem{example}[theorem]{Example}
\newtheorem{definition}[theorem]{Definition}
\newtheorem{remark}[theorem]{Remark}

\newtheorem{ques}[theorem]{Question}

\newcommand{\R}{{\mathbb{R}}}
\newcommand{\Z}{{\mathbb{Z}}}
\newcommand{\C}{{\mathbb{C}}}
\newcommand{\CP}{{\mathbb{CP}}}
\newcommand{\RP}{{\mathbb{RP}}}
\newcommand{\bS}{\mathbb{S}}

\newcommand{\gw}{\underline{c}}

\newcommand{\bigasterisk}{\mbox{\normalfont\Large\bfseries *}}

\def\eps{{\varepsilon}}

\begin{document}

\title { On the Existence of  Symplectic Barriers}
	\author{Pazit Haim-Kislev, Richard Hind, Yaron Ostrover}
	\maketitle
	\begin{abstract}
	In this note we establish the existence of a new type of rigidity of symplectic embeddings coming from obligatory intersections with symplectic planes. In particular, we prove that if a Euclidean ball is symplectically embedded in the Euclidean unit ball, then it must intersect a sufficiently fine grid of two-codimensional pairwise disjoint symplectic planes. Inspired by analogous terminology for Lagrangian submanifolds, we refer to these obstructions as  symplectic barriers.     
	\end{abstract}

\section{Introduction and Results}

Obstructions to symplectic embeddings, i.e., smooth embeddings between symplectic manifolds that respect the symplectic forms, lie at the heart of symplectic topology, and are intensively studied since the mid 1980s (see, e.g., the recent survey~\cite{schlenk-serv}). The first obstruction beyond a volume constraint is Gromov's celebrated non-squeezing theorem \cite{gromov} which states that
\begin{center} 
	\textit{	$B^{2n}(r) \overset{\mathrm s}{\hookrightarrow } Z^{2n}(1)$ if and only if $r \leq 1$,  }
\end{center}
where $B^{2n}(r)$ denotes the $2n$-dimensional Euclidean ball with radius $r$, $Z^{2n}(r)$  denotes the cylinder $B^2(r) \times \R^{2n-2}$,
and $\overset{\mathrm s}{\hookrightarrow }$ stands for a symplectic embedding.
On the other hand, a powerful measure-wise approximation result by Katok~\cite{katok} states that for every $r>0$, there are symplectic embeddings of $B^{2n}(r)$ up to sets of arbitrary small measure into $Z^{2n}(1)$.
A natural question arising from the above results, which was recently investigated in \cite{quant-nonsqueezing}, is the following
\begin{ques} \label{motivation_question}
How much does one need to remove from the ball $B^{2n}(r)$ so that it symplectically embeds into $Z^{2n}(1)$?
\end{ques}
In \cite{quant-nonsqueezing} the authors addressed this question  
when a Lagrangian plane is removed from the ball in ${\mathbb R}^4$.
In addition, they showed that if $r>1$, then the Minkowski dimension of any removed set is at least 2.  Another example of obstruction to a symplectic embedding was provided in~\cite{biran}, where it was shown that one cannot embed a ball of radius greater than or equal to $1/2$ into $\CP^n$ without intersecting $\RP^n$. In other words, the largest ball one can embed into $\CP^n \setminus \RP^n$ is of radius $1/2$.
Further examples of Lagrangian submanifolds were given in~\cite{biran} with the property that the Gromov width of their complement is smaller than the Gromov width of the ambient manifold. These were named ibid.  {\it Lagrangian barriers}. 
For some other related results see, e.g.,~\cite{Biran-Cornea, brendel-schlenk,Lee-Oh-Vianna, Traynor} and the references therein.

\medskip

The above mentioned results regarding obstructions to symplectic embeddings mostly concern the removing of a Lagrangian submanifold. 
It is natural to ask whether there are such obstructions arising from removing symplectic submanifolds.
A partial negative answer to this question can be found in \cite{mcduff-polterovich}, where it is shown that any closed complex sub-manifold in $\CP^n$ is not an obstruction to ball embeddings.  To the  best of our knowledge,  there are no known examples of embedding obstruction coming from intersections with symplectic submanifolds. In this work we show
the existence of a new type of rigidity coming from obligatory intersections with symplectic codimension two hyperplanes. In particular, this partially answers Question~\ref{motivation_question} above.

\medskip

Before we state our main result, let us recall the notion of {\it symplectic capacities}~\cite{hofer-ekeland}, which quantify obstructions to symplectic embeddings
\begin{definition}
	A  symplectic capacity is a map which associates an element of $[0, \infty]$ to every symplectic manifold $(M,\omega)$, with the following properties:
	\begin{itemize}
	\item  $c(M,\omega) \leq c(N, \tau)$ if $(M,\omega)  \overset{\mathrm s}{\hookrightarrow } (N,\tau)$ (Monotonicity)
	\item  $c(M,\alpha \omega) = |\alpha| c(M,\omega)$ for all $\alpha \in \R$, $\alpha \neq 0$ (Conformality)
	\item $c(B^{2n}(1)) = \pi = c(Z^{2n}(1)) $ (Nontriviality and Normalization)
	\end{itemize}
\end{definition}
\noindent Two examples of  symplectic capacities which naturally arise from Gromov's non-squeezing theorem
are the Gromov width and the cylindrical capacity:
\begin{align*}
\gw(M) = \sup \{\pi r^2 : B^{2n}(r) \xhookrightarrow{\mathrm s} M \}, \ 
\overline{c}(M) = \inf \{\pi r^2 : M \xhookrightarrow{s} B^{2}(r) \times \R^{2n-2} \}.
\end{align*}
Symplectic capacities can be used to formulate and generalize the questions about obstructions to symplectic embeddings in the following way:
\begin{center}
	\textit{For a submanifold $N \subset M$, what is the relation between $c(M \setminus N)$ and $c(M)$?}
\end{center}
Note that Question~\ref{motivation_question}
above can be formulated using the   cylindrical capacity. A discussion on the case where $c(M \setminus N) = c(M)$ for the special case of the Hofer-Zender capacity~\cite{hofer-zehnder} can be found in \cite{tokieda}.

\medskip

We are now in a position to state our main result:

\begin{theorem}
	\label{embed_thm}
\begin{sloppypar}	
For every $\delta>0$ and $n>1$ there exists a finite union  of codimension two 
 pairwise disjoint  symplectic  hyperplanes $\Sigma$ such that for every symplectic capacity one has $c(B^{2n}(1)\setminus\Sigma) < \pi \delta^2$. In particular, every symplectic embedding
  ${B^{2n}(\delta) \xhookrightarrow{\mathrm s} B^{2n}(1)}$ of the ball of radius $\delta$ into the unit ball must intersect $\Sigma$. 
  \end{sloppypar}
\end{theorem}

To the best of our knowledge, the rigidity coming from obligatory intersections with the symplectic hyperplanes described in Theorem~\ref{embed_thm} is a new phenomenon. Analogous to Biran's notion of Lagrangian barriers~\cite{biran}, Theorem~\ref{embed_thm} shows the existence of {\it symplectic barriers}. 
We remark that it would be interesting to better understand 
symplectic barriers, 
and in particular to further explore obstructions to symplectic embeddings coming from intersections with high-codimensional submanifolds which are not Lagrangians.

\medskip

\noindent{\bf Idea of the proof:} 
Let ${\mathbb R}^{2n} \simeq {\mathbb C}^n$ be equipped with coordinates $(z_1,\ldots,z_n)$ and the standard symplectic form $\omega = {\frac i 2} dz \wedge d \overline{z}$. 
Let $\Sigma_\varepsilon$ be the union of the symplectic codimension two hyperplanes of the form $\C^{n-1} \times \{p\} \subset \C^n$, where $p \in {\eps} \Z^{2}$. 
Fix $L > 1$. Following ideas from \cite{hind}, one can embed any convex domain $D$ minus  $\Sigma_\eps$ into a rescaling of $A^L D$, where $A^L$ is a linear map that takes $z_n$ to $L z_n$, and leaves $z_i$ fixed for $1\leq i \leq n-1$.
In other words,
there is $\beta >1$ such that
$$ D \setminus \Sigma_\eps \overset{\mathrm s}{\hookrightarrow } \beta A^L D. $$
Moreover, for a given $D$, one can show that $\beta \rightarrow 1$ when $\eps \rightarrow 0$. 
Next, one can find a symplectic linear transformation $U \in Sp(2n)$ for which  the image $D := UB^{2n}(1)$ of the ball $B^{2n}(1)$ is such that $\overline{c}(A^L D) $ is arbitrarily small. From this, one can conclude that  $\overline{c}(D \setminus \Sigma_\eps)$ is arbitrarily small for an appropriate choice of $\eps$.  Theorem~\ref{embed_thm} now follows for the pre-image of $\Sigma_\eps$ under the symplectic linear map $U$. For the exact details see Section~\ref{sec-proof-main} below.

\medskip

\begin{remark} \label{rmk:existance-of-Sigma}
Although the existence of $\Sigma$ in Theorem~\ref{embed_thm} 
is not proven by an explicit construction, one can give examples 
of such $\Sigma$'s by using specific symplectic linear transformations of the ball (see Example~\ref{ex} below). 
\end{remark}

\begin{remark}
Section 3 of the second author's paper \cite{hind} claims that in certain cases embeddings of balls can in fact be displaced from unions of symplectic hypersurfaces. However the proof given there assumes a general setting, and so our Theorem \ref{embed_thm} shows it to be incorrect. This issue will be addressed elsewhere.
\end{remark}

We note that the hypersurfaces $\Sigma$ in Theorem \ref{embed_thm} cannot be complex (with respect to the standard complex structure on ${\mathbb C}^n$). Indeed, we have the following, see also Theorem 3.1.A in \cite{mcduff-polterovich}.

\begin{proposition}
    \label{prop-complex-hyperplane}
    Let $J$ be the standard linear complex structure on $\R^{2n}$.
    For any union of $J$-holomorphic symplectic hyperplanes $\Sigma$ one has $c(B^{2n}(1) \setminus \Sigma) = c(B^{2n}(1)) $ for any symplectic capacity $c$. 
\end{proposition}

\begin{remark}
The linear symplectic group $\mathrm{Sp}(2n, \R)$ acts transitively on symplectic subspaces, and so in particular there exists a symplectic isomorphism mapping the hyperplanes $\Sigma$ from Theorem \ref{embed_thm} to complex hyperplanes. However such isomorphisms do not preserve $B^{2n}(1)$. The subgroup of $\mathrm{Sp}(2n, \R)$ preserving $B^{2n}(1)$ is the unitary group $U(n)$, which maps complex hyperplanes to themselves.
\end{remark}

\noindent {\bf Structure of the paper:}
In Section~\ref{sec-proof-main} we prove Theorem~\ref{embed_thm}. In Section~\ref{sec-example-and-proposition} we describe the example mentioned in Remark~\ref{rmk:existance-of-Sigma} above and prove Proposition~\ref{prop-complex-hyperplane}. Finally, in Section~\ref{sec-proofs-lemmas} we prove Lemma~\ref{lemma_extension} needed for the proof of Theorem~\ref{embed_thm}.

\medskip 

\noindent {\bf Acknowledgements:} We thank Leonid Polterovich for illuminating discussions and valuable comments.
P. H-K. and Y.O.  were partially supported by the ISF grant No. 938/22, and R.H. by Simons foundation grant number 663715.

\section{Proof of the Main Result} \label{sec-proof-main}
\begin{figure}[H]
\centering
	\includegraphics[width=0.8\textwidth]{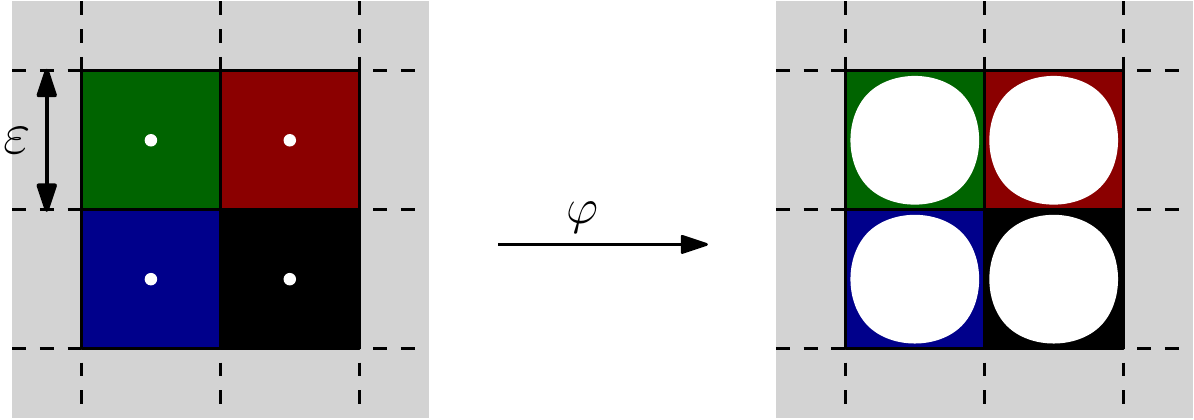}
	\caption{Illustration of the map $\varphi$.}
	\label{2-dim-map-img}
\end{figure}

\begin{figure}[h]
\centering
\includegraphics[width=0.55\textwidth, trim = 2cm 5.7cm 1cm 5cm]{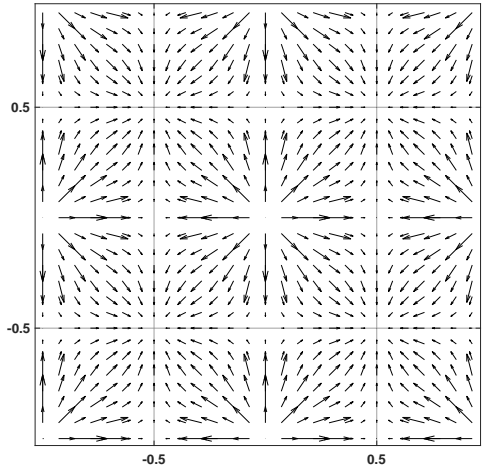}
 \caption{An illustration of the vector field constructed in the proof of Lemma \ref{lemma_extension}.}
 \label{vf_fig}
 \end{figure}

Before we prove Theorem \ref{embed_thm}, let us lay the groundwork.
Let $\eps > 0$, and denote by $\{ G_{\alpha} \} \subset \R^2$ the collection of all squares with sides of length $\eps$ and centers in $\eps \Z^2$.
We shall need the following lemma, the proof of which is  given in Section~\ref{sec-proofs-lemmas} (see also Figure~\ref{2-dim-map-img}).
\begin{lemma}
\label{lemma_extension}
Let $L \geq 1$. For every $\varepsilon > 0$ there is an embedding $$ \varphi : {\mathbb R}^2 \setminus \varepsilon {\mathbb Z}^2 \to  {\mathbb R}^2 \setminus \varepsilon   {\mathbb Z}^2, $$ 
which preserves the grid squares $G_{\alpha}$ and is 
such that
$L^2 \varphi^* \omega = \omega$.
\end{lemma}
The idea of the proof of Lemma~\ref{lemma_extension} is to construct (by hand) a vector field that contracts $\omega$, pushing points near each puncture radially outwards. The condition on the grid squares, and integrability of the vector field, follow by ensuring that the vector field is tangent to the boundary of the $G_{\alpha}$, and vanishes at their vertices (see Figure \ref{vf_fig}).

\medskip 

Next, denote by $\Sigma_\eps$ the union of the symplectic codimension two hyperplanes in ${\mathbb R}^{2n} \simeq \C^n$, 
\begin{align} \label{def-std-grid}
\Sigma_\eps := \bigcup \{(z_1,z_2,\ldots,z_n) \in \C^n : z_n \in \eps \Z^{2} \}.
\end{align}

\begin{proof}[{\bf Proof of Theorem \ref{embed_thm}}]
Note that from the monotonicity property of symplectic capacities it is enough to prove the theorem for the cylindrical capacity. 
For a convex body $K \subset \R^{2n}$ define  $$\lambda(K) := \inf_{u \in \bS^{2n-1}} h_K(u), \ {\rm and} \ h_K(u) := \sup \{ \langle x, u \rangle : x \in K \}.$$ Note that $h_K$ is the support function of $K$, and for centrally symmetric convex bodies, i.e., when $K=-K$, it measures the width of $K$ in a given direction.
Let \begin{equation} \label{al}
 A^L : {\mathbb R}^{2n} \rightarrow {\mathbb R}^{2n}, \, (z_1, \dots, z_n) \mapsto (z_1, \dots, z_{n-1}, Lz_n).  
\end{equation}

\vskip 5pt

First we wish to show that for any convex 
$D \subset \R^{2n}$ containing the origin,
and for every $\eps>0$ and $L>1$, one can find a symplectic embedding 

\begin{equation} \label{required-embedding-Psi}
   \Psi: D \setminus \Sigma_\eps \xhookrightarrow{\mathrm s} \left(1+\frac{\sqrt{2}\eps L}{\lambda(A^L D)}\right) A^L D.
\end{equation}
 In what follows, we will use the above fact for $D$ being a specific ellipsoid.
To establish the embedding $\Psi$, denote by $\Pi_{i_1,\ldots,i_j}$ the projection to the complex coordinates $i_1,\ldots,i_j$.
Consider the following approximation of the set $D$ (see Figure~\ref{d_eps_fig}):
$$ D_\eps := \bigcup_\alpha \left(  \{\Pi_{1,\ldots,n-1}(x) : x \in D, \, \Pi_{n}(x) \in G_\alpha\}  \times G_\alpha \right), $$
where the $G_{\alpha} \subset \R^2$ 
is the collection of all squares with sides of length $\eps$ and centers in $\eps \Z^2$.

\begin{figure}[h]
\centering
	\includegraphics[width=0.5\textwidth]{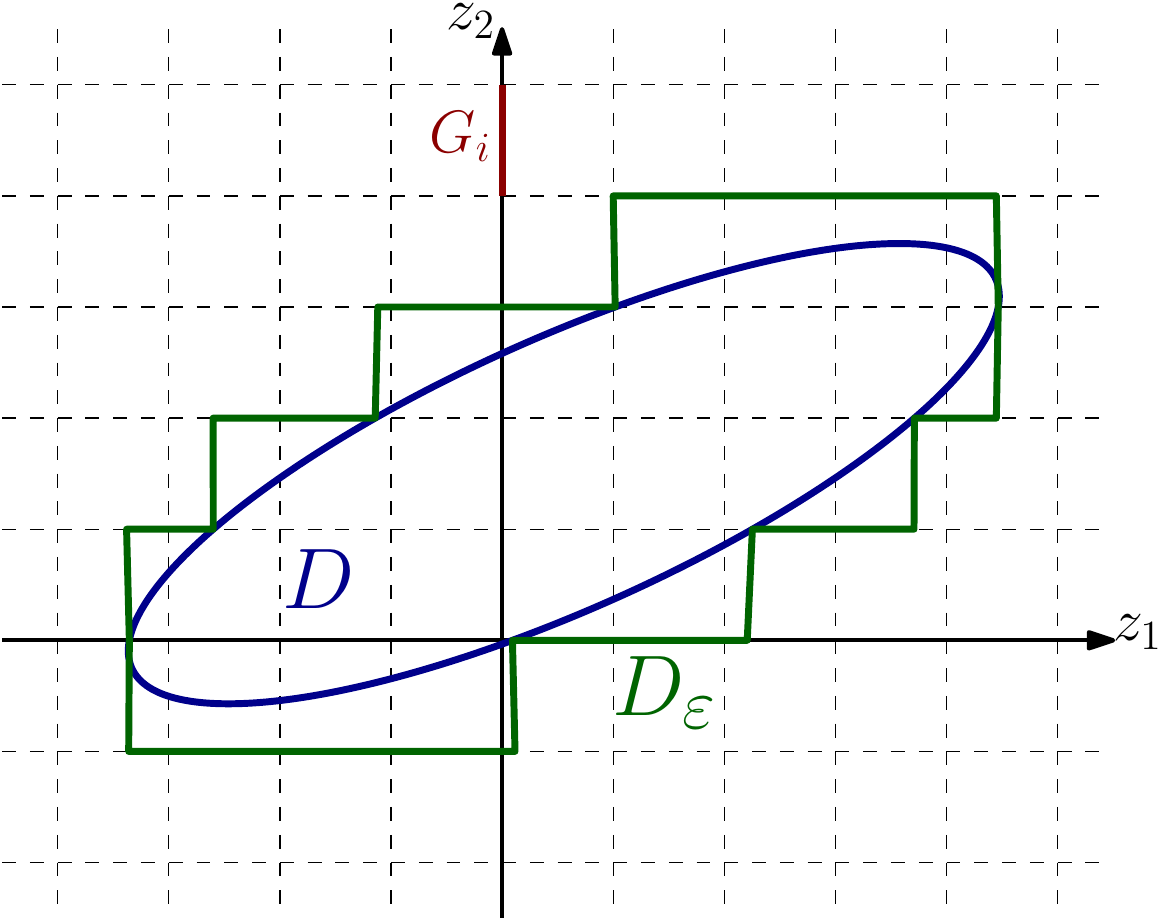}
 \caption{A 2-dimensional example for the approximated set $D_\eps$.}
 \label{d_eps_fig}
 \end{figure}
Note that since $d_{\text{H}}(A^L D_\eps, A^L D) \leq \sqrt{2}\eps L$, where $d_{\text{H}}$ is the Hausdorff distance, one has
$$ A^L D_\eps \subseteq A^L D + B^{2n}(\sqrt{2}\eps L) \subseteq \left(1+\frac{\sqrt{2}\eps L}{\lambda(A^L D)}\right) A^L D.$$
Hence, it is enough to prove that one can symplectically embed $D \setminus \Sigma_\epsilon$ into $A^L D_\eps$.
To  this end, let $\phi: \R^2 \setminus \eps\Z^2 \to \R^2 \setminus \eps \Z^2 $ be the map from Lemma \ref{lemma_extension}. Then $\psi = L\phi$ is a symplectic embedding $\psi: \R^2 \setminus \eps\Z^2 \to \R^2 \setminus L\eps \Z^2 $, mapping the squares $G_{\alpha}$ to $LG_{\alpha}$. Hence $\Psi = {\rm Id} \times \cdots \times {\rm Id} \times \psi$ is also a symplectic map, with the property that

\begin{equation} \label{the-Psi-embedding}
  \Psi(D \setminus \Sigma_\eps) \subset \Psi(D_\eps \setminus \Sigma_\eps) \subset A^L D_\eps \subset 
\left(1+\frac{\sqrt{2}\eps L}{\lambda(A^L D)}\right) A^L D  
\end{equation}
as required in~\eqref{required-embedding-Psi}.

\medskip

Next, using a lemma by Eliashberg~\cite{eliashberg} (cf. Lemma 2 in~\textsection 2.2 in~\cite{HZ-book}),  one has that for every $L>1$ and every $a > 0$, there are symplectic matrices $U$ and $V$ such that
\[
V A^L U=\left(\begin{array}{c|c}
\begin{array}{cc}
a & 0\\
0 & a
\end{array} & 0\\
\hline \\*
\bigasterisk & \,\,\,\,\,\bigasterisk\,\,\,\,\,\\
\\
\end{array}\right).
\]
We remark that the above mentioned lemma actually holds more generally, replacing $A^L$ by any non symplectically conformal linear isomorphism of ${\mathbb R}^{2n}$; we do not use this generalization.
From Eliashberg's lemma we conclude that $V A^L U B^{2n}(1) \subset Z^{2n}(a) $, where $U$ and $V$ are linear symplectic maps. This implies that for $D := U B^{2n}(1)$ one has
$$ \overline{c}(A^L D) = \overline{c}(V A^L U B^{2n}(1)) \leq \pi a^2. $$

To complete the proof of Theorem~\ref{embed_thm}, let $\delta>0$. Fix $L>1$, and choose $a < \delta$ and $\eps$ small enough to get
$$ \overline{c}(D \setminus \Sigma_\eps) \leq \left( 1 + \frac{\sqrt{2} \eps L}{\lambda(A^L D)} \right)^2 \overline{c}(A^L D) \leq \left( 1 + \frac{\sqrt{2} \eps L}{\lambda(A^L D)} \right)^2 \pi a^2  < \pi \delta^2,
$$
where the first inequality follows from the existence of $\Psi$ and monotonicity of symplectic capacities, the second inequality from our choice of $D$, and the third by taking $\epsilon$ sufficiently small.
Hence
 $\Sigma := U^{-1} \Sigma_\eps$ is a union of codimension two symplectic hyperplanes with
$$ \overline{c}(B(1) \setminus \Sigma) = \overline{c}( D \setminus \Sigma_\eps) < \pi \delta^2$$
as required. The proof of Theorem~\ref{embed_thm} is now complete. 
 \end{proof}

\section{Symplectic vs. Complex Hyperplanes} \label{sec-example-and-proposition}
In this section we attempt to better understand the symplectic barriers $\Sigma$ from Theorem~\ref{embed_thm}. We start with explicit examples for such $\Sigma$'s, which are symplectic linear images of the standard grid~\eqref{def-std-grid}, 
and continue with the proof of Proposition~\ref{prop-complex-hyperplane} which shows that $\Sigma$ cannot be taken to be complex.

\begin{example} \label{ex}
We shall use the notations from the proof of Theorem~\ref{embed_thm} above. For simplicity, we assume in what follows that $n=2$. The arguments can be easily extended to higher dimensions.

\medskip 

We wish to find a symplectic matrix 
$S $
such that $B^{4}(1) \setminus S^{-1} \Sigma_\varepsilon$ has smaller capacity then $B^4(1)$, where $\Sigma_{\varepsilon}$ is the standard grid~\eqref{def-std-grid}. 
From the proof of Theorem \ref{embed_thm}, in particular by~(\ref{the-Psi-embedding}), one has that for every $L>1$
	\begin{equation} \label{main-eq-in-example}  
 c( B^4(1) \setminus S^{-1} \Sigma_\eps) \leq \left( 1 + \frac{\sqrt{2} \eps L}{\lambda(A^L S B^4(1))} \right)^2 c(A^L  S B^4(1)), \end{equation}
 for any symplectic capacity $c$, where $A^L$ is the linear transformation given by \eqref{al}.
Since $\varepsilon>0$ can be taken to be arbitrarily small, in what follows we focus on bounding the capacity $c(A^L S B^4(1))$. 

\medskip 

We start with an informal discussion 
describing the motivation for the selection of the specific symplectic matrix $S$ below. We follow the guiding principle that, for very large $L$, the ellipsoid $A^L S B^4(1)$ will have small capacity when the areas of the intersections of the  ellipsoid $S B^4(1)$ with the planes $\{z_2 = b\}$, $b \in \C$, are sufficiently small.
This follows roughly speaking from the fact that, when $L$ is large, a linear symplectic transformation mapping the ellipsoid $A^L S B^4(1)$ to a canonical form $\{ \frac{|z_1|^2}{a} + \frac{|z_2|^2}{b} <1\}$ will approximately preserve the axis $\{z_2=0\}$, and so $a$, and hence the capacity, must also be small.
We remark that an ellipsoid with fixed capacity and arbitrarily small intersections was suggested by Polterovich (see Proposition 8.2.4 in~\cite{schlenk} and Question 3.6 in~\cite{mcduff}). More precisely, one has  
$$ \inf_{S \in {\rm Sp(4, \R)}} \sup_{b \in \C} {\rm Area}(S(B^4(1)) \cap \{z_2 = b\}) = 0 .$$
In fact, for $\alpha >0$, consider the symplectic matrix 
	$$M_\alpha :=
	\left(
	\begin{matrix}
	\alpha & 0 & 0 & -1 \\
	0 & 1/\alpha & 0 & 0 \\
	0 & -1 & \alpha & 0 \\
	0 & 0 & 0 & 1/\alpha
	\end{matrix}
	\right). $$
A direct computation shows that $$\sup_{b \in \C} {\rm Area}(M_\alpha B^4(1) \cap \{z_2 = b\}) = {\rm Area}(M_\alpha B^4(1) \cap \{z_2 = 0\}) = \frac{\pi \alpha}{\alpha^2 + 1}.$$
This implies in particular that every compact domain in $\C^2$ has a linear symplectic image  whose intersections with the symplectic planes $\C \times \{pt\}$ have arbitrarily small area.
 
 \medskip
Based on the above informal discussion we  take $S = M_{\alpha}$, where $\alpha = 0.5$, set $L=2$, and estimate the upper bound in~$(\ref{main-eq-in-example})$. Note that the capacity of the ellipsoid $A^L S B^4(1)$ can be computed based on the symplectic normal form for ellipsoids (see Section 1.7 in~\cite{HZ-book}), and one has that 
$c(A^L S B^4(1)) \approx \pi(0.7423)^2$.
Moreover, a direct computation gives that $\lambda(A^L S B^4(1)) \approx 0.4849$. Hence, from~$(\ref{main-eq-in-example})$ it follows that 
$$ c( B^4(1) \setminus S^{-1} \Sigma_\eps) \leq \left( 1 + \frac{\sqrt{2} \eps L}{\lambda(A^L SB^4(1))} \right)^2 c(A^L S B^4(1)) \approx \pi (1+5.8335 \eps)^2 (0.7423)^2 ,$$
and for $\eps$ small enough we get a specific set of symplectic hyperplanes which are barriers as required.

\medskip
	
\end{example}

\begin{proof}[{\bf Proof of Proposition~\ref{prop-complex-hyperplane}}]
Note that for the proof it is enough to find a symplectic embedding of the ball $B^{2n}(r)$ into $B^{2n}(1) \setminus \Sigma$ for any $r <1$.
 This can be done by constructing, for every $r<1$, a compactly supported isotopy  of $\Sigma$ such that the image of $\Sigma$ does not intersect the ball $B^{2n}(r)$, and moreover, throughout the homotopy $\Sigma$ remains symplectic. Indeed, from a Moser's type argument and Banyaga's isotopy extension theorem (see Theorem 3.3.2 in~\cite{Mc-Sal}) it follows that the above mentioned isotopy can be extended to a symplectomorphism of the ball. 
 
\medskip

 To construct the required isotopy, we first assume, by replacing $B^{2n}(1)$ with a ball of radius $1-\eps$ centered at a point $z$ with $|z|<\eps$ that $\Sigma$ is disjoint from $0$.  
 Now consider the rescaling vector field $Y = \sum x_i \partial x_i + y_i \partial y_i $ with a radial cutoff $\chi(|z|)$ near the boundary of the ball. We choose $\chi(|z|)$ such that $\chi(|z|)=1$ when $|z|<r$ and $\chi(1) =0$.
	Denote by $G$ the image of the flow generated by $\chi(|z|)Y$ at some time $t$. Then for $t$ sufficiently large we have $G(\Sigma)$ disjoint from $B^{2n}(r)$. To check the symplectic condition, as $\Sigma$ is $J$-holomorphic it is enough to prove that $\omega(G_* u, G_* J u) > 0$ for every $u \neq 0$.
	Assume $u$ is a tangent vector at a point $x$. One can write $u = \lambda_1 x + \lambda_2 J x + v$ where $v \perp x$ and $v \perp J x$.
	Notice that there exist $\alpha>0$, $\beta>0$ (depending upon $x$ and $t$) such that  $G_* x = \alpha x$ and $G_*w = \beta w$ for every $w \perp x$. Then we get
	\begin{align*}
	\omega(G_* u, G_* J u) &= \omega(G_* (\lambda_1 x + \lambda_2 Jx + v), G_*(\lambda_1 J x - \lambda_2 x + Jv) ) \\
	&= \omega (\lambda_1 \alpha x + \lambda_2 \beta Jx + \beta v, \lambda_1 \beta Jx - \lambda_2 \alpha x + \beta J v ) \\
	&= \lambda_1^2 \alpha \beta \omega(x, Jx) -\lambda_2^2 \alpha \beta \omega(Jx , x)  + \beta^2 \omega(v, Jv) \\
	&= (\lambda_1^2 + \lambda_2^2)\alpha \beta \omega(x,Jx) + \beta^2 \omega(v, Jv) >0 
	\end{align*}
 This completes the proof of the proposition. 
\end{proof}

\section{Proof of the Technical Lemma}
\label{sec-proofs-lemmas}

\begin{proof}
[{\bf Proof of Lemma \ref{lemma_extension}}]
Let $\omega$ be the standard symplectic form on ${\mathbb R}^2$. 
By rescaling, it suffices to show that there exists a diffeomorphism $$\psi : \R^2 \setminus \Z^2 \to \R^2 \setminus \Z^2$$ which preserves the grid squares $Q_i$ with vertices at the points $(m+\frac 12, n + \frac 12)$, $m, n \in \Z$, and satisfies $L\psi^* \omega = \omega$ for any $L \ge 1$. To do this we look for a  vector field $X$ on $\R^2 \setminus \Z^2$, tangent to the edges of the grid squares $Q_i$, vanishing at their vertices, and with ${\cal L}_X \omega = -\omega$. 
This is equivalent to finding a $1$-form $\lambda$ with $d \lambda = - \omega$, such that the restriction of $\lambda$ to the edges of the grid squares vanishes.

Next, let $f(x)$ be a 1-periodic function with $f(n + \frac 12)=0$ and $f'(x)=\frac 12$ when $x$ is close to $n + \frac 12$, for all integers $n$. Also, assume $|f'(x)|$ is small elsewhere.
Consider the $1$-form on $\R^2$ given by $$\lambda_1 = f(y) dx - f(x) dy.$$
Then $d \lambda_1 = -(f'(x) + f'(y))dx \wedge dy$ which, by our hypotheses on $f$, is symplectic near the grid $G = \{ x + \frac 12 \in \Z \, \mathrm{or} \, y + \frac 12 \in \Z \}$, and  is minus the standard symplectic form near the corners. 
Note that with respect to $\omega$, the form $\lambda_1$ is dual to the vector field $$X_1 = -f(x) \partial_x - f(y) \partial_y.$$
We note that $X_1$ is tangent to $G$, vanishing at the vertices, and the corresponding flow
is contracting towards $G$. 
It is not hard to check that one can construct an isotopy fixing $G$, such that the pull back of $\lambda_1$ with respect to this isotopy, denoted by $\lambda_2$, satisfies that $d \lambda_2 = - \omega$
 in a neighborhood of the grid.

\vskip 5pt

Next, near each integer point, that is, near the center of the grid squares, we work in polar coordinates with center $r=0$ and define $$\mu = -(\frac 12 r^2 - \frac{1}{2\pi}) d \theta.$$
Note that $d \mu = -\omega$, and the contracting vector field $Y = (-\frac 12 r +\frac{1}{2 \pi r}) \partial_r$ points away from the lattice points $\Z^2$. We denote by $\lambda_3$ a $1$-form on $\R^2 \setminus \Z^2$ coinciding with $\lambda_2$ near $G$, and with $\mu$ near $\Z^2$. Hence, $\omega + d \lambda_3$ is a closed $2$-form on $\R^2 \setminus \Z^2$ which vanishes near $G \cup \Z^2$, and thus extends to $\R^2$. Moreover, as $\lambda_3$ vanishes when restricted to $G$, by Stokes' Theorem we see that the integral over each grid square $Q_i$ is $$\int_{Q_i}( \omega +  d \lambda_3) = 1 - \int_{\gamma} \mu =0,$$
where $\gamma$ is an infinitesimal circle around the center, oriented negatively.

\vskip 5pt
Let $\zeta$ be a primitive of $\omega + d\lambda_3$. By Poincar\'{e}'s Lemma, adding differentials of functions we may assume $\zeta$ vanishes near $\Z^2$. Also, as
$$\int_{\partial Q_i} \zeta = \int_{Q_i}( \omega +  d \lambda_3) = 0$$ for all grid squares $Q_i$, we may also add a differential so that $\zeta$ vanishes near $G$.

Finally we define $\lambda = \lambda_3 - \zeta$. This is a primitive of $-\omega$ on $\R^2 \setminus \Z^2$ such that the dual vector field has flow preserving grid squares and moving away from $\Z^2$. Hence the $\omega$ contracting flow exists for all time as required.

\end{proof}

\bibliography{references}

\begin{thebibliography}{10}

\bibitem{biran}
{\sc P.~Biran}, {\em Lagrangian barriers and symplectic embeddings}, Geom.
  Funct. Anal., 11 (2001), pp.~407--464.

\bibitem{Biran-Cornea}
{\sc P.~Biran and O.~Cornea}, {\em Rigidity and uniruling for {L}agrangian
  submanifolds}, Geom. Topol., 13 (2009), pp.~2881--2989.

\bibitem{brendel-schlenk}
{\sc J.~Brendel and F.~Schlenk}, {\em Pinwheels as {L}agrangian barriers},
  arXiv:2210.00280.

\bibitem{hofer-ekeland}
{\sc I.~Ekeland and H.~Hofer}, {\em Symplectic topology and {H}amiltonian
  dynamics}, Math. Z., 200 (1989), pp.~355--378.

\bibitem{eliashberg}
{\sc Y.~M. Eliashberg}, {\em A theorem on the structure of wave fronts and its
  application in symplectic topology}, Funktsional. Anal. i Prilozhen., 21
  (1987), pp.~65--72, 96.

\bibitem{gromov}
{\sc M.~Gromov}, {\em Pseudo holomorphic curves in symplectic manifolds},
  Invent. Math., 82 (1985), pp.~307--347.

\bibitem{hind}
{\sc R.~Hind}, {\em {Symplectic capacities of domains in {$\Bbb C^2$}}},
  International Mathematics Research Notices, 2006 (2006).
\newblock ID 37171.

\bibitem{hofer-zehnder}
{\sc H.~Hofer and E.~Zehnder}, {\em A new capacity for symplectic manifolds},
  in Analysis, et cetera, Academic Press, Boston, MA, 1990, pp.~405--427.

\bibitem{HZ-book}
{\sc H.~Hofer and E.~Zehnder}, {\em Symplectic invariants and {H}amiltonian
  dynamics}, Birkh\"{a}user Advanced Texts: Basler Lehrb\"{u}cher.,
  Birkh\"{a}user Verlag, Basel, 1994.

\bibitem{katok}
{\sc A.~B. Katok}, {\em Ergodic perturbations of degenerate integrable
  {H}amiltonian systems}, Izv. Akad. Nauk SSSR Ser. Mat., 37 (1973),
  pp.~539--576.

\bibitem{Lee-Oh-Vianna}
{\sc W.~Lee, Y.-G. Oh, and R.~Vianna}, {\em Asymptotic behavior of exotic
  {L}agrangian tori {$T_{a,b,c}$} in {$\Bbb C{\rm P}^2$} as {$a + b + c\to
  \infty$}}, J. Symplectic Geom., 19 (2021), pp.~607--634.

\bibitem{mcduff}
{\sc D.~McDuff}, {\em Fibrations in symplectic topology}, in Proceedings of the
  {I}nternational {C}ongress of {M}athematicians, {V}ol. {I} ({B}erlin, 1998),
  no.~Extra Vol. I, 1998, pp.~339--357.

\bibitem{mcduff-polterovich}
{\sc D.~McDuff and L.~Polterovich}, {\em Symplectic packings and algebraic
  geometry}, Invent. Math., 115 (1994), pp.~405--434.
\newblock With an appendix by Yael Karshon.

\bibitem{Mc-Sal}
{\sc D.~McDuff and D.~Salamon}, {\em Introduction to symplectic topology},
  Oxford Graduate Texts in Mathematics, Oxford University Press, Oxford,
  third~ed., 2017.

\bibitem{quant-nonsqueezing}
{\sc K.~Sackel, A.~Song, U.~Varolgunes, and J.~Zhu~J.}, {\em On certain
  quantifications of {G}romov's non-squeezing theorem}, To appear in Geometry
  and Topology,  (2022).
\newblock With an appendix by Jo\'{e} Brendel.

\bibitem{schlenk}
{\sc F.~Schlenk}, {\em Embedding problems in symplectic geometry}, vol.~40 of
  De Gruyter Expositions in Mathematics, Walter de Gruyter GmbH \& Co. KG,
  Berlin, 2005.

\bibitem{schlenk-serv}
\leavevmode\vrule height 2pt depth -1.6pt width 23pt, {\em Symplectic embedding
  problems, old and new}, Bull. Amer. Math. Soc. (N.S.), 55 (2018),
  pp.~139--182.

\bibitem{tokieda}
{\sc T.~F. Tokieda}, {\em Isotropic isotopy and symplectic null sets}, Proc.
  Nat. Acad. Sci. U.S.A., 94 (1997), pp.~13407--13408.

\bibitem{Traynor}
{\sc L.~Traynor}, {\em Symplectic packing constructions}, J. Differential
  Geom., 42 (1995), pp.~411--429.

\end{thebibliography}
\bibliographystyle{siam}

\vskip10pt

\noindent Pazit Haim-Kislev \\
\noindent School of Mathematical Sciences, Tel Aviv University, Israel \\
\noindent e-mail: pazithaim@mail.tau.ac.il
\vskip 10pt

\noindent Richard  Hind \\
\noindent Department of Mathematics, University of Notre Dame, IN, USA. \\
\noindent e-mail: hind.1@nd.edu
\vskip 10pt

\noindent Yaron Ostrover \\
\noindent School of Mathematical Sciences, Tel Aviv University, Israel \\
\noindent e-mail: ostrover@tauex.tau.ac.il

\end{document}